\documentclass[11pt]{amsart}
\usepackage{amsmath,amsfonts,amscd,amssymb,epsf, euscript, eufrak}
\usepackage{amsxtra,mathrsfs} 
\newtheorem{theorem}{Theorem}[section]
\newtheorem{proposition}[theorem]{Proposition}
\newtheorem{lemma}[theorem]{Lemma}
\newtheorem{corollary}[theorem]{Corollary}
\theoremstyle{definition}
\newtheorem{definition}[theorem]{Definition}

\theoremstyle{remark} \newtheorem{remark}[theorem]{Remark}
\numberwithin{equation}{section}

\newcommand{\complex}[1]{\mathsf{#1}} 

\newcommand{\SSS}{\complex{S}}

\newcommand{\KKK}{\complex{K}}

 \DeclareMathOperator{\Diff}{Diff}
\DeclareMathOperator{\Mob}{M\ddot{o}b}

 \DeclareMathOperator{\PSL}{PSL}
\DeclareMathOperator{\PSU}{PSU}

\newcommand{\Del}{\mathbb{D}}

\newcommand{\U}{{\mathbb{U}}}

\newcommand{\R}{{\mathbb{R}}}
\newcommand{\C}{{\mathbb{C}}}
\newcommand{\bL}{\mathbb{L}}

\newcommand{\curly}[1]{\mathscr{#1}}

\newcommand{\cB}{\curly{B}}

\newcommand{\cP}{\curly{P}}

\newcommand{\bk}{\backslash}
\newcommand{\pa}{\partial}

\newcommand{\ov}{\overline}

\newcommand{\vep}{\varepsilon}
\newcommand{\f}{\mathrm{f}}
\newcommand{\g}{\mathrm{g}}
\newcommand{\z}{\bar{z}}

\begin{document}
\title[Weil-Petersson potential on
the Quasi-Fuchsian deformation space]{A different expression of
the Weil-Petersson potential on the Quasi-Fuchsian deformation
space}

\author{ Lee-Peng Teo}\address{Faculty of Information Technology \\
Multimedia University \\ Jalan Multimedia, Cyberjaya\\
63100, Selangor\\Malaysia} \subjclass[2000]{ 32G15, 30F60, 30F10}
\keywords{Potential, Teichm\"uller space, quasi-Fuchsian}
\email{lpteo@mmu.edu.my}

\begin{abstract}We extend a definition of the Weil-Petersson potential on the
universal Teichm\"uller space to the quasi-Fuchsian deformation
space. We prove that up to a constant, this function coincides
with the Weil-Petersson potential on the quasi-Fuchsian
deformation space. As a result, we prove a lower bound for the
potential on the quasi-Fuchsian deformation space.\end{abstract}
\date{\today}
\maketitle

\section{Introduction} In
\cite{TT2}, we defined a Hilbert manifold structure on the
universal Teichm\"uller space $T(1)$. Under this structure, $T(1)$
is a disjoint union of uncountably many components. We denoted by
$T_0(1)$ the component that contains the identity element. It can
 be characterized as the completion of the space
$\Mob(S^1)\bk\Diff_+(S^1)$ under the Weil-Petersson metric. Hence
it is the largest submanifold of $T(1)$ where the Weil-Petersson
metric can be defined. The Weil-Petersson metric on
$\Mob(S^1)\bk\Diff_+(S^1)$, introduced by Kirillov \cite{Ki} via
the orbit method, has been of interest to both mathematicians and
physicists. It is a right-invariant K\"ahler metric, and hence it
may play some role in the canonical quantization of the space
$\Mob(S^1)\bk\Diff_+(S^1)$.

In our subsequent work \cite{TT3}, we defined a Weil-Petersson
potential on $T_0(1)$ in two different ways, and showed that they
are equal. The first definition of the Weil-Petersson potential is
given by  $\SSS_1: T_0(1) \rightarrow \R$,
\begin{align*} \SSS_1([\mu]) =
\iint\limits_{\Del}\left|\frac{(\f^{\mu})^{\prime\prime}}{(\f^{\mu})^{\prime}}\right|^2d^2z+
\iint\limits_{\Del^*}
\left|\frac{\g_{\mu}^{\prime\prime}}{\g_{\mu}^{\prime}}\right|^2
d^2z -4\pi\log|\g_{\mu}'(\infty)|,
\end{align*}
where $w_{\mu}=\g_{\mu}^{-1}\circ\f^{\mu}$ is the conformal
welding corresponding to $[\mu]\in T_{0}(1)$. The second
definition of the Weil-Petersson potential comes from the study of
the Grunsky operator $K_1$ of the univalent function $\f^{\mu}$
associated to a point $[\mu]\in T(1)$. We proved that the Grunsky
operator associated to $[\mu]$ is Hilbert-Schmidt if and only if
$[\mu]\in T_0(1)$. Hence the function $\SSS_2 : T_0(1)\rightarrow
\R$ given by $$\SSS_2([\mu]) = \log \det (I- K_1K_1^*)$$ is
well-defined. We proved that $$\SSS_2 =-\frac{1}{12 \pi}\SSS_1.$$

We call $\SSS_1$ the universal Liouville action in our paper
\cite{TT3}. In fact, given $\Gamma$ a cocompact quasi-Fuchsian
group, the quasi-Fuchsian deformation space $\mathfrak{D}(\Gamma)$
can be canonically mapped into the universal Teichm\"uller space,
with the totally real submanifold
--- the space of Fuchsian groups --- mapped to the point identity. The
function $\SSS_1$ bears a lot of resemblance to the critical value
of the Liuoville action functional $S_{cl}$ we constructed in
\cite{TT1}. However, we haven't established the precise relation
between them. In this paper, we are going to extend the definition
of the function $\SSS_2$ to the quasi-Fuchsian deformation space
of a quasi-Fuchsian group $\Gamma$, and prove that up to
constants, it coincides with the critical Liouville action
$S_{cl}$. As a result, we show that (see Corollary
\ref{inequality})
$$S_{cl}([\mu]) \leq 8\pi (2g-2),$$ with equality appears if and
only if $[\mu]$ corresponds to a Fuchsian group.

It is our intention to keep this paper concise. Hence we will not
repeat the background material and conventions. We refer them to
our previous papers \cite{TT1, TT2, TT3}.

\section{Definition of the function $S_2$ on quasi-Fuchsian deformation
space}\label{DS2}

\subsection{The quasi-Fuchsian deformation space $\mathfrak{D}_g$}
\subsubsection{ A model of $\mathfrak{D}_g$}
 We fix a model for the quasi-Fuchsian
deformation space of genus $g\geq 2$ in the following way. Let
$\Gamma\in \PSU(1,1)$ be a normalized cocompact Fuchsian group of
genus $g$. Let $\mathcal{A}^{-1,1}(\Gamma)$ be the space of
bounded Beltrami differentials for $\Gamma$ and
$\mathcal{B}^{-1,1}(\Gamma)$ the unit ball of
$\mathcal{A}^{-1,1}(\Gamma)$ with respect to the sup-norm. For
each $\mu \in \mathcal{B}^{-1,1}(\Gamma)$, there exists a unique
quasi-conformal (q.c.) mapping $w_{\mu}: \hat{\C}\rightarrow
\hat{\C}$ satisfying the Beltrami equation
$$(w_{\mu})_{\z}=\mu (w_{\mu})_z$$ and fixing the points $-1,-i,
1$. The quasi-Fuchsian deformation space of the Fuchsian group
$\Gamma$ is defined as
$$\mathfrak{D}(\Gamma) =
\mathcal{B}^{-1,1}(\Gamma)/\sim,$$ where $\mu \sim \nu$ if and
only if $w_{\mu}\bigr|_{S^1}=w_{\nu}\bigr|_{S^1}$.

 Given $[\mu]\in \mathcal{B}^{-1,1}(\Gamma)$, let $\Gamma^{\mu}
 = w_{\mu} \circ \Gamma\circ w_{\mu}^{-1}$. By definition,
 $\Gamma^{\mu}$ is
 a normalized quasi-Fuchsian group and it
 is a Fuchsian group if and only if
 $\mu$ is symmetric,
 i.e.
$$\mu(z) = \ov{\mu\left(\frac{1}{\z}\right)}\frac{z^2}{\z^2}.$$
 There is a canonical isomorphism
$\mathfrak{D}(\Gamma) \xrightarrow{\simeq}
\mathfrak{D}(\Gamma^{\mu})$ given by $[\nu] \mapsto [\lambda]$,
where $\lambda$ is the Beltrami differential of $w_{\nu}\circ
w_{\mu}^{-1}$. We define $\mathfrak{D}_g$, the quasi-Fuchsian
deformation space of genus $g$ to be the space
$\mathfrak{D}(\Gamma)$, where $\Gamma$ is a Fuchsian group of
genus $g$\footnote{We can let $\Gamma$ to be any quasi-Fuchsian
group of genus $g$. But for the description of some properties of
$\mathfrak{D}(\Gamma)$ in terms of equivalence classes of Beltrami
differentials, it will be convenient to assume that $\Gamma$ is a
Fuchsian group.}. It is unique up to isomorphism.

Given $[\mu]\in \mathfrak{D}(\Gamma)$, let $\Omega_1=
w_{\mu}(\Del)$ and $\Omega_2=w_{\mu}(\Del^*)$. Then
$\Omega=\Omega_1\sqcup \Omega_2$ is the domain of discontinuity of
$\Gamma^{\mu}$ acting on $\hat{\C}$ and $\mathcal{C} =
w_{\mu}(S^1)$ is the quasi-circle separating the domains
$\Omega_1$ and $\Omega_2$. There exists a unique q.c. mapping $f:
\hat{\C}\rightarrow \hat{\C}$ such that $f|_{\Del}$ is
holomorphic, $f$ fixes the points $-1,-i,1$  and
$f(\Del)=\Omega_1$. Similarly, there exists a unique q.c. mapping
$g: \hat{\C}\rightarrow \hat{\C}$ such that $g|_{\Del^*}$ is
holomorphic, $g$ fixes the points $-1,-i,1$  and
$g(\Del^*)=\Omega_2$. By abusing notation, we also denote by $f$
and $g$ the univalent functions $f\bigr|_{\Del}$ and
$g\bigr|_{\Del^*}$. We say that $(f,g)$ is the pair of univalent
functions associated to the point $[\mu]\in \mathfrak{D}_g$. In
fact, $(f,g)$ is independent of the choice of the model
$\mathfrak{D}(\Gamma)$ for $\mathfrak{D}_g$. It only depends on
the quasi-Fuchsian group $\Gamma^{\mu}$. In case $\Gamma^{\mu}$ is
a Fuchsian group, $\Omega_1=\Del$, $\Omega_2=\Del^*$ and hence
$f=\text{id}$ and $g=\text{id}$. Using the biholomorphisms $f$ and
$g$, we define the pair of Fuchsian groups $(\Gamma_1, \Gamma_2)$
associated to $\Gamma^{\mu}$ by $\Gamma_1=f^{-1}\circ
\Gamma^{\mu}\circ f$ and $\Gamma_2=g^{-1}\circ \Gamma^{\mu}\circ
g$.

Given $[\nu]\in \mathfrak{D}_g$, let $\mu_1, \ldots,\mu_d$ be a
basis of $\Omega^{-1,1}(\Gamma^{\nu})$. The Bers coordinates at
the point $[\nu]$ is defined by the correspondence
$\boldsymbol{\vep}=(\vep_1, \ldots, \vep_d) \mapsto
\Gamma^{\boldsymbol{\vep}}= w_{\boldsymbol{\vep}}\circ
\Gamma^{\nu}\circ w_{\boldsymbol{\vep}}^{-1}$, where
$w_{\boldsymbol{\vep}}:\hat{\C}\rightarrow \hat{\C}$ is the unique
q.c. mapping with Beltrami differential
$\vep_1\mu_1+\ldots\vep_d\mu_d$ and fixing the points $-1,-i,1$.

\subsubsection{The tangent and cotangent space of
$\mathfrak{D}_g$}

The holomorphic tangent space at the point $[\nu]\in
\mathfrak{D}_g$ is isomorphic to the space of harmonic Beltrami
differentials $\Omega^{-1,1}( \Gamma^{\nu})$ of $\Gamma^{\nu}$.
 The holomorphic
cotangent space at the point $[\mu]\in \mathfrak{D}_g$ is
isomorphic to the vector space $\Omega^{2,0}(\Gamma^{\mu})$ of
holomorphic quadratic differentials of $\Gamma^{\mu}$.

Given $\mu \in \Omega^{-1,1}(\Gamma^{\nu})$, we denote by
$\frac{\pa}{\pa\vep_{\mu}}$ and $\frac{\pa}{\pa\bar{\vep}_{\mu}}$
the holomorphic and anti-holomorphic vector fields in a
neighbourhood of $[\nu]$ defined using the Bers coordinates at the
point $[\nu]$.

\subsubsection{The embedding $T_g \hookrightarrow \mathfrak{D}_g$
and the map $\mathfrak{D}_g\rightarrow T(1)$} Let $T_g$ be the
Teichm\"uller space of genus $g$ and $T(1)$ the universal
Teichm\"uller space. $T_g$ can be realized as a complex
submanifold of the quasi-Fuchsian deformation space
$\mathfrak{D}_g =\mathfrak{D}(\Gamma)$. Namely
$$ T_g=T(\Gamma)=\left\{ [\mu]\in \mathfrak{D}(\Gamma) \;:\; \mu(z)=0
\;\text{on}\; \Del\cup S^1.\right\}.$$ We denote by $i:T_g
\hookrightarrow \mathfrak{D}_g$ the canonical inclusion map.

To define the map $\Xi: \mathfrak{D}(\Gamma)\rightarrow T(1)$, we
use the following model for $T(1)$.
$$T(1) = L^{\infty}(\Del^*)_1/\sim,$$ where given $\mu\in
L^{\infty}(\Del^*)_1$, a Beltrami differential on $\Del^*$ with
sup-norm less than $1$, we extend $\mu$ to be zero outside
$\Del^*$ and let $w^{\mu}: \hat{\C}\rightarrow \hat{\C}$ to be the
unique q.c. mapping with Beltrami differential $\mu$ and fixing
the points $-1,-i,1$. $\mu\sim \nu$ if and only if
$w^{\mu}\bigr|_{\Del}=w^{\nu}\bigr|_{\Del}$.  Given $[\mu]\in
\mathfrak{D}_g$, let $(f,g)$ be the pair of univalent functions
associated to $[\mu]$. Extend $f$ and $g$ to be q.c. mappings. The
map $\Xi: \mathfrak{D}(\Gamma)\rightarrow T(1)$ is defined so that
$\Xi([\mu])$ is the equivalence class of the Beltrami differential
of $g\bigr|_{\Del^*}$. It is independent of the choice of the
representative $\mu$ of  $[\mu]$. It is not a one-to-one mapping.
In particular, the subspace of Fuchsian groups in
$\mathfrak{D}(\Gamma)$ is mapped to the point $[0]\in T(1)$. If
$(f,g)$ is the pair of univalent functions  associated to
$[\mu]\in \mathfrak{D}_g$, there exists a linear fractional
transformation $\lambda\in\PSL(2,\C)$ such that the functions $\f=
\lambda\circ f$ and  $\g=\lambda\circ g$ satisfies
$$\f(0)=0,
\hspace{1cm}\f'(0)=1,\hspace{1cm}\g(\infty)=\infty.$$The pair of
functions $(\f, \g)$
 is
then the functions in the conformal welding $w=\g^{-1}\circ \f $
associated to $\Xi([\mu])\in T(1)$.

  Obviously, we also have a canonical complex analytic embedding
  $\tilde{i}: T_g \hookrightarrow T(1)$. It factors as $\tilde{i}=
  \Xi\circ i$.

\subsubsection{The inversion on $\mathfrak{D}_g$}

There is an inversion map $\mathfrak{I}: \mathfrak{D}_g\rightarrow
\mathfrak{D}_g$ on the quasi-Fuchsian deformation space induced by
the inversion $\iota: \hat{\C}\rightarrow \hat{\C}$, $z\mapsto
1/\bar{z}$. It is defined in the obvious way:
$$\mathfrak{I}([\mu]) = [\iota^* \mu], \hspace{1cm} \text{for}\; [\mu]\in \mathfrak{D}_g,  $$
where $$ \iota^*\mu(z) =
\ov{\mu\left(\frac{1}{\z}\right)}\frac{z^2}{\z^2}.$$ The space of
Fuchsian groups is the set of fixed points of this
map\footnote{This is only true when we use a Fuchsian group
$\Gamma$ for the model $\mathfrak{D}_g=\mathfrak{D}(\Gamma)$.}. By
uniqueness of q.c. mappings, it is easy to see that
\begin{align}\label{inverse} w_{\iota^* \mu}
(z)=\frac{1}{\ov{w_{\mu}(1/\z)}},\hspace{1cm} f[\iota^*\mu](z) =
\frac{1}{\ov{g[\mu](1/\z)}}, \hspace{1cm} g[\iota^*\mu](z) =
\frac{1}{\ov{f[\mu](1/\z)}}.\end{align}The relations between the
quasi-Fuchsian groups and associated Fuchsian groups are given by
$$\Gamma^{\iota^*\mu}=\iota\circ \Gamma^{\mu}\circ \iota \hspace{1cm} \Gamma_1[\iota^*\mu]
=\iota\circ \Gamma_2[\mu]\circ \iota,
\hspace{1cm}\Gamma_2[\iota^*\mu]=\iota\circ \Gamma_1[\mu]\circ
\iota.$$

\subsection{The function $S_2$}

\subsubsection{Integral operators associated to $(f,g)$}
Let
\begin{align*}
A^1_{2}(\Del)&= \left\{\psi \;\text{holomorphic on}\,\, \Del: \,
\|\psi\|_{2}^2 = \iint\limits_{ \Del}
\left|\psi(z)\right|^2d^2z < \infty \right\},\\
A^1_{2}(\Del^*)&= \left\{\psi \;\text{holomorphic on}\,\, \Del^*:
\, \|\psi\|_{2}^2 = \iint\limits_{ \Del^*}
\left|\psi(z)\right|^2d^2z < \infty \right\}
\end{align*}
be Hilbert spaces of holomorphic functions on $\Del$ and $\Del^*$
and denote by $\ov{A_2^1(\Del)}$ and $\ov{A_2^1(\Del^*)}$ the
corresponding Hilbert-spaces of anti-holomorphic functions. Given
a pair $f:\Del\rightarrow \C$, $g:\Del^*\rightarrow \hat{\C}$ of
univalent functions such that $\hat{\C}\setminus(f(\Del)\cup
g(\Del^*))$ has measure zero , define the kernel functions
\begin{align*}
K_{1}(z,w) &= \frac{1}{\pi}\left(\frac{1}{(z-w)^2}-
\frac{f'(z)f'(w)}{(f(z)-f(w))^2}\right),\\
K_2(z,w)&= \frac{1}{\pi}\frac{f'(z)
g'(w)}{(f(z)-g(w))^2},\\
K_{3}(z,w) &= \frac{1}{\pi}
\frac{g'(z)f'(w)}{(g(z)-f(w))^2},\\
K_4(z,w)&=\frac{1}{\pi}\left( \frac{1}{(z-w)^2}-\frac{g'(z)
g'(w)}{(g(z)-g(w))^2}\right).
\end{align*}
They define linear operators $K_l$, $l=1,2,3,4$ as follows,
\begin{alignat*}{3}
K_1 &: \ov{A^{1}_{2}(\Del)} \rightarrow A^{1}_{2}(\Del), &\qquad
(K_1 \psi)(z) &= \iint\limits_{\Del}
K_1(z,w)\ov{\psi(w)}d^2w,\\
K_2&: \ov{A^{1}_{2}(\Del^*)}\rightarrow A^{1}_{2}(\Del), &\qquad
(K_2\psi)(z)&= \iint\limits_{\Del^*}
K_2(z,w)\ov{\psi(w)}d^2w,\\
K_3&: \ov{A^{1}_{2}(\Del)} \rightarrow A^{1}_{2}(\Del^*), &\qquad
(K_3\psi)(z)&= \iint\limits_{\Del}
K_3(z,w)\ov{\psi(w)}d^2w,\\
K_4&: \ov{A^{1}_{2}(\Del^*)} \rightarrow A^{1}_{2}(\Del^*),
&\qquad (K_4\psi)(z)&= \iint\limits_{\Del^*}
K_4(z,w)\ov{\psi(w)}d^2w.
\end{alignat*}
The generalized Grunsky equality says that these operators satisfy
the following relations (see e.g., \cite{TT3}):
\begin{align}\label{relation2}
K_1 K_1^* + K_2 K_2^* = I,\hspace{1cm} K_3 K_1^* + K_4 K_2^*=0,\\
K_1 K_3^* + K_2 K_4^* = 0,\hspace{1cm}  K_3 K_3^* +
K_4K_4^*=I.\nonumber
\end{align}
$K_2, K_3$ are invertible operators and $K_1, K_4$ are operators
of norm strictly less than one.

\begin{remark}
Our definition of the operators $K_l$ here can be viewed as the
'pull--back' of the corresponding definition on $T(1)$ via the map
$\Xi:\mathfrak{D}_g\rightarrow T(1)$.
\end{remark}

\begin{remark}\label{singular} If $\phi$ is a
holomorphic function on $\Del$, the principal--valued integral
\begin{align*}
\iint\limits_{\Del} \frac{\ov{\phi(w)}}{(z-w)^2}d^{2}w
\end{align*}
vanishes identically. Hence we can also represent operators $K_1$
and $K_4$ by the singular kernels
\begin{align*}
-\frac{1}{\pi}\frac{f'(z)f'(w)}{(f(z)-f(w))^2}\quad\text{and}\quad
-\frac{1}{\pi}\frac{g'(z)g'(w)}{(g(z)-g(w))^2}.
\end{align*}
\end{remark}
\subsubsection{The definition of $S_2$.} Let $\KKK_l=K_lK_l^*$, $l=1,2,3,4$. In \cite{TT3},
we define the function $\SSS_2:T_0(1)\rightarrow \R$, which up to
a multiplicative constant is a Weil-Petersson potential on
$T_0(1)$, by
$$\SSS_2=\log\det(I-\KKK_1)=\log\det(I-\KKK_4).$$ This definition
cannot be extended to $\mathfrak{D}_g$ since if $\KKK_1$ and
$\KKK_4$ are defined using the pair $(f,g)$ associated to a point
on $\mathfrak{D}_g$, they are not trace-class operators unless
$\KKK_1=\KKK_4=0$ (see the proof in \cite{TT3}), which corresponds
to the case $\Gamma^{\mu}$ is a Fuchsian group. Nevertheless,
motivated by the series expansion
$$\log\det(I-K)= - \text{Tr}\left(\sum_{n=1}^{\infty}\frac{K^n}{n}\right)$$
valid when the operator $K$ has norm less than $1$, we want to
consider the following operators
\begin{align*}\mathfrak{O}_1 =\sum_{n=1}^{\infty}\frac{\KKK_1^n}{n}
\hspace{1cm}\text{and}\hspace{1cm}\mathfrak{O}_2
=\sum_{n=1}^{\infty}\frac{\KKK_4^n}{n}.\end{align*}

\begin{lemma}
The operators $\mathfrak{O}_1$ and $\mathfrak{O}_2$ are
well-defined operators with kernels
\begin{align*}
\mathfrak{O}_1(z,w) = \sum_{n=1}^{\infty}
\frac{\mathcal{K}_{1,n}(z,w)}{n}
\hspace{1cm}\text{and}\hspace{1cm} \mathfrak{O}_2(z,w) =
\sum_{n=1}^{\infty} \frac{\mathcal{K}_{4,n}(z,w)}{n},
\end{align*}
which converge absolutely and uniformly on compact subsets of
$\Del\times \Del$ and $\Del^*\times \Del^*$ respectively. Here
$\mathcal{K}_{i,n}(z,w)$, $i=1,4$ is the kernel of the operator
$\KKK_i^n$.

\end{lemma}
\begin{proof}
It is sufficient to consider the operator $\mathfrak{O}_1$. First,
we notice that for $n\geq 2$,
\begin{align*}
\mathcal{K}_{1,n}(z,w) =&
\iint\limits_{\Del}\iint\limits_{\Del}K_1(z,\zeta)\Bigl(\left(K_1^*K_1\right)
^{n-1}\Bigr)(\zeta,\eta)
K_1^*(\eta, w) d^2\zeta d^2\eta \\
=&\iint\limits_{\Del}\iint\limits_{\Del}
\Bigl(\left(K_1K_1^*\right)^{n-1}\Bigr)(\eta,\zeta)K_1(\zeta,z)
\ov{K_1(\eta,w)}d^2\zeta d^2\eta\\=&\left\langle \KKK_1^{n-1}v_z,
v_w\right\rangle.
\end{align*}
Here we denote by $\langle \; .\;, \;.\;\rangle$ the inner product
on the Hilbert space $A_2^1(\Del)$, and $v_z$ is the holomorphic
function
$$v_z(\zeta)=K_1(z,\zeta)$$with norm
$$\Vert v_z\Vert_2^2=\KKK_1(z,z).$$By Cauchy-Schwarz inequality,
we have
\begin{align*}
\left|\mathcal{K}_{1,n}(z,w) \right| \leq \Vert
\KKK_1^{n-1}v_z\Vert_2\Vert v_w\Vert_2\leq \Vert
\KKK_1\Vert_{\infty}^{n-1}\Vert v_z\Vert_2\Vert v_w\Vert_2.
\end{align*}
Hence
\begin{align*}
 \sum_{n=1}^{\infty}\left|\frac{\mathcal{K}_{1,n}(z,w)}{n}\right|\leq
 \left(
 \sum_{n=1}^{\infty} \frac{\left\Vert
 \KKK_1\right\Vert_{\infty}^{n-1}}{n}\right) \Vert v_z\Vert_2\Vert v_w\Vert_2
\end{align*}
which converges absolutely and uniformly on compact subsets of
$\Del\times \Del$ since $\Vert \KKK_1\Vert_{\infty}<1$ and
$$ \Vert v_z\Vert_2^2 \leq \frac{1}{\pi(1-|z|^2)^2}.$$
(see \cite{TT3}). In fact, as a function of $w$,
$$\Vert \mathcal{K}_{1,n}(z,w)\Vert_2^2= \mathcal{K}_{1, 2n}(z,z)
\leq \Vert \KKK_1\Vert_{\infty}^{2n-1} \Vert v_z\Vert_2^2,$$ which
implies that as a function of $w$,
$$\lim_{k\rightarrow \infty}\sum_{n=1}^{k}
\frac{\mathcal{K}_{1,n}(z,w)}{n}=\sum_{n=1}^{\infty}
\frac{\mathcal{K}_{1,n}(z,w)}{n}$$ in $A_2^1(\Del)$.
 Hence it follows that if $\psi\in A_2^1(\Del)$,
\begin{align*}
\left(\mathfrak{O}_1\bar{\psi}\right)(z) =&\lim_{k\rightarrow
\infty}\sum_{n=1}^k
\frac{(\KKK_1^n\bar{\psi})(z)}{n}=\lim_{k\rightarrow
\infty}\iint\limits_{\Del} \left(\sum_{n=1}^{k}
\frac{\mathcal{K}_{1,n}(z,w)}{n}\right)\ov{\psi(w)}d^2w \\
&=\iint\limits_{\Del} \left(\sum_{n=1}^{\infty}
\frac{\mathcal{K}_{1,n}(z,w)}{n}\right)\ov{\psi(w)}d^2w,
\end{align*}
which proves that $\sum_{n=1}^{\infty}
\frac{1}{n}\mathcal{K}_{1,n}(z,w)$ is the kernel for
$\mathfrak{O}_1$.

\end{proof}

\begin{corollary}\label{pos}
Let $(f,g)$ and $(\Gamma_1, \Gamma_2)$ be the pairs of univalent
functions and Fuchsian groups associated to a point on
$\mathfrak{D}_g$. The functions $\mathfrak{O}_1(z,z)$ and
$\mathfrak{O}_2(z,z)$ defined using $(f,g)$ are nonnegative real
valued continuous functions on $\Del$ and $\Del^*$ that are
automorphic $(1,1)$ forms with respect to $\Gamma_1$ and
$\Gamma_2$ respectively.

\end{corollary}

\begin{proof}
Again, it suffices to consider $\mathfrak{O}_1$.  It follows from
the proof of the previous lemma that
$$\mathfrak{O}_1(z,z) = \sum_{n=1}^{\infty}
\frac{\mathcal{K}_{1,n}(z,z)}{n}$$ converges absolutely and
uniformly on compact subsets of $\Del$. Hence it is continuous.
Moreover, since $\KKK_1$ is a positive self-adjoint operator,
$$
\mathcal{K}_{1,n}(z,z) =\left\langle \KKK_1^{n-1}v_z,
v_z\right\rangle\geq 0$$ for all $n$. Hence
$\mathfrak{O}_1(z,z)\geq 0$. Now for every $\gamma\in \Gamma_1$,
there exists $\tilde{\gamma}\in \Gamma^{\mu}$ such that
$$f\circ\gamma= \tilde{\gamma}\circ f.$$ Hence it is easy to
check from the definition of $K_1(z,w)$ that
$$ K_1(\gamma z, \gamma w)\gamma'(z)\gamma'(w)= K_1(z,w)\hspace{1cm}\forall \gamma\in \Gamma_1$$
and
consequently
$$ \mathcal{K}_{1,n}(\gamma z, \gamma w) \gamma'(z)\ov{\gamma'(w)}
=\mathcal{K}_{1,n}(z,w)\hspace{1cm}\forall \gamma\in \Gamma_1.$$
Therefore, $\mathfrak{O}_1(z,z)$ is an automorphic $(1,1)$-form
with respect to $\Gamma_1$.
\end{proof}

Now we can define the function $S_2: \mathfrak{D}_g \rightarrow
\R$.
\begin{definition}
The function $S_2: \mathfrak{D}_g \rightarrow \R$ is defined as
follows:
\begin{displaymath}
S_2([\mu]) = \iint\limits_{\Gamma_1\bk\Del}
\mathfrak{O}_1(z,z)d^2z=\sum_{n=1}^{\infty}
\frac{1}{n}\iint\limits_{\Gamma_1\bk\Del}
\mathcal{K}_{1,n}(z,z)d^2z ,
\end{displaymath}
where $\mathfrak{O}_1(z,w)$ is defined using the univalent
function $f$ associated to $[\mu]$.
\end{definition}

\begin{remark}
The function $S_2$ can be considered as the regularized trace of
the operator $-\log\det(I-\KKK_1)$ on $A_2^1(\Del)$.
\end{remark}

\subsubsection{Behavior of $S_2$ under inversion}
 The relations \eqref{inverse} give us the following relations for the
 kernels $K_l$ associated to $[\mu]$ and $[\iota^*\mu]$ on
 $\mathfrak{D}_g$:
$$K_l[\iota^*\mu](z,w)=\ov{K_{4-l}[\mu]\left(\frac{1}{\z},
\frac{1}{\bar{w}}\right)}\frac{1}{z^2}\frac{1}{w^2},
\hspace{1cm}l=1,2,3,4.$$ In particular,
$$\mathfrak{O}_4[\mu](z,z)=\mathfrak{O}_1[\iota^*\mu]\left(\frac{1}{\z},\frac{1}{\z}\right)\frac{1}{|z|^4}
.$$ Using this, we have the following result.
\begin{proposition}\label{invstable}
The function $S_2$ is invariant under inversion, i.e. $S_2\circ
\mathfrak{I} =S_2$.
\end{proposition}
\begin{proof}
Given a point $[\mu]$ on $\mathfrak{D}_g$ with the associated
univalent functions $(f,g)$ and Fuchsian groups
$(\Gamma_1,\Gamma_2)$, we are going to prove that
$$\iint\limits_{\Gamma_1 \bk\Del}\mathcal{K}_{1,n}(z,z) d^2z
=\iint\limits_{\Gamma_2\bk\Del^*} \mathcal{K}_{4,n}(z,z)d^2z$$ for
all $n$.  From the relations \eqref{relation2} we
have\begin{align}\label{op1} \KKK_1(z,w) +
\KKK_2(z,w)&=I_1(z,w)=\frac{1}{\pi(1-z\bar{w})^2},
\hspace{1cm}z,w\in \Del;\\
\KKK_3(z,w) + \KKK_4(z,w)&=I_2(z,w)=\frac{1}{\pi(1-z\bar{w})^2},
\hspace{1cm}z,w\in \Del^*;\nonumber
\end{align}
where $I_1$ and $I_2$ are the identity operators on $A_2^1(\Del)$
and $A_2^1(\Del^*)$ respectively. On the other hand,
\begin{align*}
\mathcal{K}_{1,n}(z,w) =& \Bigl(\left(
I_1-\KKK_{2}\right)^{n}\Bigr)(z,w)=\sum_{k=0}^n(-1)^k
\begin{pmatrix} n\\k\end{pmatrix} \left(\KKK_{2}^k\right)(z,w), \\
\mathcal{K}_{4,n}(z,w) =& \Bigl(\left(
I_2-\KKK_{3}\right)^{n}\Bigr)(z,w)=\sum_{k=0}^n(-1)^k
\begin{pmatrix} n\\k\end{pmatrix} \left(\KKK_{3}^k\right)(z,w).
\end{align*}
Now for $k\geq 1$,
\begin{align*}
&\iint\limits_{\Gamma_1\bk\Del}\left(\KKK_{2}^k\right)(z,z)d^2z
=\iint\limits_{\Gamma_1\bk\Del}\;\;\iint\limits_{\Del^*}
K_2(z,\zeta)\left((K_2^*K_2)^{k-1}K_2^*\right)(\zeta, z)d^2\zeta
d^2z\\
=&\sum_{\gamma_2\in \Gamma_2}\;\;
\iint\limits_{\Gamma_1\bk\Del}\;\;\iint\limits_{\Gamma_2\bk\Del^*}
K_2(z,\gamma_2\zeta)\left((K_2^*K_2)^{k-1}K_2^*\right)(\gamma_2\zeta,
z)|\gamma_2'(\zeta)|^2 d^2\zeta d^2z.
\end{align*}
For every $\gamma_2\in \Gamma_2$, there exists $\gamma\in
\Gamma^{\mu}$ such that
$$ \gamma\circ g= g\circ \gamma_2$$ and $\gamma_1\in \Gamma_1$ such
that
$$\gamma\circ f= f\circ\gamma_1.$$
Hence
$$K_2(\gamma_1 z, \gamma_2 \zeta)\gamma_1'(z)\gamma_2'(\zeta) =
K_2(z,\zeta) $$ whenever the pair of elements $\gamma_1\in
\Gamma_1,\gamma_2\in \Gamma_2$  are associated to the same element
$\gamma\in \Gamma$. Then it is easy to show that
$$\left((K_2^*K_2)^{k-1}K_2^*\right)(\gamma_2\zeta,
z)\ov{\gamma_2'(\zeta)} =\left((K_2^*K_2)^{k-1}K_2^*\right)(\zeta,
\gamma_1^{-1}z)\ov{(\gamma_1^{-1})'(z)}.$$ Consequently,
\begin{align*}
&\iint\limits_{\Gamma_1\bk\Del}\left(\KKK_{2}^k\right)(z,z)d^2z\\
&\sum_{\gamma_1\in \Gamma_1}\;\;
\iint\limits_{\Gamma_1\bk\Del}\;\;\iint\limits_{\Gamma_2\bk\Del^*}
K_2(\gamma^{-1}_1
z,\zeta)\left((K_2^*K_2)^{k-1}K_2^*\right)(\zeta,
\gamma_1^{-1}z)\left|(\gamma_1^{-1})'(z)\right|^2 d^2\zeta d^2z\\
=&\iint\limits_{\Gamma_2\bk\Del^*}
\iint\limits_{\Del}K_2(z,\zeta)\left((K_2^*K_2)^{k-1}K_2^*\right)(\zeta,
z)d^2z d^2\zeta \\
=&\iint\limits_{\Gamma_2\bk\Del^*} \iint\limits_{\Del}K_3(\zeta,
z)\left((K_3^*K_3)^{k-1}K_3^*\right)(z, \zeta)d^2z
d^2\zeta=\iint\limits_{\Gamma_2\bk\Del^*}\left(\KKK_3^k\right)(z,z)d^2z.
\end{align*}
We have used the equality $K_2(z,w)= K_3(w,z) \forall z\in \Del,
w\in \Del^*$ in the last line. Finally, from \eqref{op1} we have
\begin{align*}
\iint\limits_{\Gamma_1\bk\Del} \mathcal{K}_{1,n}(z,z)d^2z =&
\frac{1}{4\pi}\text{Area}(\Gamma_1\bk \Del)
+\sum_{k=1}^{n}(-1)^k\begin{pmatrix} n\\k\end{pmatrix}
\iint\limits_{\Gamma_1\bk\Del}\left(\KKK_{2}^k\right)(z,z)d^2z\\
=&\frac{1}{4\pi}\text{Area}(\Gamma_2\bk \Del^*)
+\sum_{k=1}^{n}(-1)^k\begin{pmatrix} n\\k\end{pmatrix}
\iint\limits_{\Gamma_2\bk\Del^*}\left(\KKK_{3}^k\right)(z,z)d^2z\\=&\iint\limits_{\Gamma_2\bk\Del^*}
\mathcal{K}_{4,n}(z,z)d^2z.
\end{align*}
It follows from the definition that
$$ \iint\limits_{\Gamma_1[\mu]\bk\Del} \mathfrak{O}_1[\mu](z,z)d^2z=
\iint\limits_{\Gamma_2[\mu]\bk\Del^*}
\mathfrak{O}_4[\mu](z,z)d^2z=\iint\limits_{\Gamma_1[\iota^*\mu]\bk\Del}
\mathfrak{O}_1[\iota^*\mu](z,z)d^2z.$$

\end{proof}

\section{The first variation of the function $S_2$}

Given $\mu\in\Omega^{-1,1}(\Gamma^{\nu})$ a tangent vector at the
point $[\nu]$, we define
$$\mu_1 = f^*(\mu|_{\Omega_1}) \hspace{1cm}\text{and}\hspace{1cm}
\mu_2=g^*(\mu|_{\Omega_2}).$$

We separate the computation of the variation of $S_2$ into a few
lemmas.
\begin{lemma}\label{lemma2}Given $[\nu]\in \mathfrak{D}_g$, let
$\mu\in \Omega^{-1,1}(\Gamma^{\nu})$ be such that $\mu$ has
support on $\Omega_2$. Let $(f^{\vep}, g_{\vep})$ be the univalent
functions associated to $\Gamma^{\vep}=w_{\vep\mu}\circ
\Gamma^{\nu}\circ w_{\vep\mu}^{-1}$. At the point $[\nu]$, the
variation of the kernel $\KKK_1$ in the direction $\mu$ is given
by
\begin{align*}
&\frac{\pa}{\pa \vep}\Bigr|_{\vep=0} \KKK_1^{\vep}(z,w)
=-\iint\limits_{\Del^*}\iint\limits_{\Del^*}\mu_2(u)K_2(z,u)K_4(u,\zeta)
K_2^*(\zeta, w)d^2ud^2\zeta
\end{align*}
\end{lemma}
\begin{proof}See the proof of Lemma 2.7 and Theorem 3.1 in
\cite{TT3}.
\end{proof}

\begin{lemma}\label{lemmaite}
Given $[\nu]\in \mathfrak{D}_g$, let $\mu\in
\Omega^{-1,1}(\Gamma^{\nu})$ be such that $\mu$ has support on
$\Omega_2$. Let $(f^{\vep}, g_{\vep})$ and
$(\Gamma_1^{\vep},\Gamma_2^{\vep})$ be the univalent functions and
Fuchsian groups associated to $\Gamma^{\vep}=w_{\vep\mu}\circ
\Gamma^{\nu}\circ w_{\vep\mu}^{-1}$. Then for all $n\geq 1$,
\begin{align*}
\frac{\pa}{\pa \vep}\Bigr|_{\vep=0}
\iint\limits_{\Gamma_1^{\vep}\bk \Del}\mathcal{K}_{1,n}^{\vep}
(z,z) d^2z=n\iint\limits_{\Gamma_1\bk
\Del}\iint\limits_{\Del}\left(\frac{\pa}{\pa
\vep}\Bigr|_{\vep=0}\KKK_1^{\vep}(z, \zeta)\right) \mathcal{K}_{1,
n-1}(\zeta,z)d^2\zeta d^2z.
\end{align*}

\end{lemma}
\begin{proof}
Since $f^{\vep}=w_{\vep}\circ f$, the group $\Gamma_1^{\vep}=
(f^{\vep})^{-1}\circ \Gamma^{\vep} \circ f^{\vep}$ is a constant
with respect to $\vep$, i.e.
$\Gamma_1^{\vep}=\Gamma_1^{0}=\Gamma_1$. From Lemma \ref{lemma2},
it is easy to check that for all $\gamma_1\in \Gamma_1$,
\begin{align*}
\left(\frac{\pa}{\pa \vep}\Bigr|_{\vep=0}
\KKK_1^{\vep}\right)(\gamma_1 z,\gamma_1 w)\gamma_1'(z)
\gamma_1'(w) = \frac{\pa}{\pa \vep}\Bigr|_{\vep=0}
\KKK_1^{\vep}(z,w).
\end{align*}
We find that
\begin{align*}
&\frac{\pa}{\pa \vep}\Bigr|_{\vep=0}
\iint\limits_{\Gamma_1^{\vep}\bk
\Del}\mathcal{K}_{1,n}^{\vep} (z,z) d^2z\\
=& \frac{\pa}{\pa \vep}\Bigr|_{\vep=0}\iint\limits_{\Gamma_1\bk
\Del}\iint\limits_{\Del}\dots\iint\limits_{\Del} \KKK_1^{\vep}(z,
\zeta_1)\KKK_1^{\vep}(\zeta_1, \zeta_2)\ldots
\KKK_1^{\vep}(\zeta_{n-1}, z) d^2\zeta_1\ldots
d^2\zeta_{n-1}d^2z\\
=&\iint\limits_{\Gamma_1\bk
\Del}\iint\limits_{\Del}\dots\iint\limits_{\Del}\left(\frac{\pa}{\pa
\vep}\Bigr|_{\vep=0}\KKK_1^{\vep}(z,
\zeta_1)\right)\KKK_1(\zeta_1, \zeta_2)\ldots \KKK_1(\zeta_{n-1},
z) d^2\zeta_1\ldots
d^2\zeta_{n-1}d^2z\\
&+\ldots\\
&+\iint\limits_{\Gamma_1\bk
\Del}\iint\limits_{\Del}\dots\iint\limits_{\Del}\KKK_1(z,
\zeta_1)\KKK_1(\zeta_1, \zeta_2)\ldots \left(\frac{\pa}{\pa
\vep}\Bigr|_{\vep=0}\KKK_1^{\vep}(\zeta_{n-1}, z)\right)
d^2\zeta_1\ldots d^2\zeta_{n-1}d^2z\\
=&n\iint\limits_{\Gamma_1\bk
\Del}\iint\limits_{\Del}\left(\frac{\pa}{\pa
\vep}\Bigr|_{\vep=0}\KKK_1^{\vep}(z, \zeta)\right) \mathcal{K}_{1,
n-1}(\zeta,z)d^2\zeta d^2z.
\end{align*}
\end{proof}

\begin{lemma}\label{lemma4}
Given $[\nu]\in \mathfrak{D}_g$, let $\mu\in
\Omega^{-1,1}(\Gamma^{\nu})$ be such that $\mu$ has support on
$\Omega_2$. Let $(f^{\vep}, g_{\vep})$ and
$(\Gamma_1^{\vep},\Gamma_2^{\vep})$ be the univalent functions and
Fuchsian groups associated to $\Gamma^{\vep}=w_{\vep\mu}\circ
\Gamma^{\nu}\circ w_{\vep\mu}^{-1}$. Then there exists an $r>0$
such that the series
$$\sum_{n=1}^{\infty}\frac{1}{n}\frac{\pa}{\pa\vep}\iint\limits
_{\Gamma_1^{\vep}\bk\Del} \mathcal{K}^{\vep}_{1,n}(z,z)d^2z,$$
converges uniformly to
$$ -\iint\limits_{
\Del^*}\iint\limits_{\Gamma_1\bk\Del}\iint\limits_{\Del}\iint\limits_{\Del^*}
\mu^{\vep}_2(u)K^{\vep}_2(z, u)K^{\vep}_4(u,
\eta)(K_2^{\vep})^*(\eta,\zeta)(I-\KKK_1^{\vep})^{-1}(\zeta,z)d^2\eta
d^2\zeta d^2z d^2u$$ in the ball $\{\vep\in \C\;:\; |\vep|<r\}.$
Here $\mu_2^{\vep}$ is the Beltrami differential
$$g_{\vep}^*\left(\left(\frac{\mu}{1-|\vep\mu|^2}\frac{(w_{\vep\mu})_z}{\ov
{(w_{\vep\mu})}_{\z}}\right)\circ (w_{\vep\mu})^{-1}\right).$$
\end{lemma}
\begin{proof}
By shifting the origin of differentiation, it is easy to see from
Lemma \ref{lemma2} and Lemma \ref{lemmaite} that
\begin{align*}
&\frac{\pa}{\pa\vep}\iint\limits_{\Gamma_1\bk\Del}
\mathcal{K}^{\vep}_{1,n}(z,z)d^2z\\=&-n\iint\limits_{\Gamma_1\bk
\Del}\iint\limits_{\Del}\iint\limits_{\Del^*}\iint\limits_{\Del^*}\mu_2^{\vep}
(u)K_2^{\vep}(z,u)K_4^{\vep}(u,\eta)(K_2^{\vep})^*(\eta,\zeta)\mathcal{K}_{1,
n-1}^{\vep}(\zeta,z)d^2u d^2\eta  d^2\zeta d^2z.
\end{align*}
Now using the property of the operators $K_l$, $l=1,2,3,4$, we
find that the $\ell^2$--norm of the function
$$u_{w}^{\vep}(\zeta) = \iint\limits_{ \Del^*}
\iint\limits_{\Del^*} \mu_2^{\vep} (u)K_2^{\vep}(w,u)K_4^{\vep}(u,
\eta)(K_2^{\vep})^*(\eta, \zeta)d^2ud^2\eta $$ satisfies the
inequality
\begin{align*}
\Vert u_{w}^{\vep}\Vert_2^2 \leq &\iint\limits_{\Del^*}\left|
\iint\limits_{\Del^*}
\mu_2^{\vep}(u) K_2^{\vep}(w,u)K_4^{\vep}(u,\zeta)d^2u\right|^2 d^2\zeta\\
\leq & \iint\limits_{\Del^*} \left|\mu_2^{\vep}(\zeta)\right|^2
|K_2^{\vep}(w, \zeta)|^2 d^2\zeta\\
\leq & \Vert \mu_2^{\vep}\Vert_{\infty}^2 \KKK_2^{\vep}(w,w)\leq
\frac{\Vert
\mu_2^{\vep}\Vert_{\infty}^2}{\pi(1-|w|^2)^2}.\end{align*} Hence,
the $\ell^2$--norm of the function
\begin{align*}
v_{w,n}^{\vep}(z)=&\iint\limits_{\Del}\iint\limits_{\Del^*}\iint\limits_{\Del^*}\mu_2^{\vep}
(u)K_2^{\vep}(w,u)K_4^{\vep}(u,\eta)(K_2^{\vep})^*(\eta,\zeta)d^2u
d^2\eta \mathcal{K}_{1, n-1}^{\vep}(\zeta,z)d^2\zeta\\
=& \bigl(((K_1^*K_1)^{\vep})^{n-1} u_{w}^{\vep}\bigr)(z)
\end{align*}
satisfies
$$ \Vert v_{w,n}^{\vep}\Vert_2 \leq \Vert
(K_1^*K_1)^{\vep}\Vert_{\infty}^{n-1} \Vert
u_{w,n}^{\vep}\Vert_2\leq \frac{\Vert
K_1^{\vep}\Vert_{\infty}^{2n-2}\Vert
\mu_2^{\vep}\Vert_{\infty}}{\sqrt{\pi}(1-|w|^2)}.
$$
Therefore, by Lemma 2.3 in \cite{TT3}, we have
\begin{align*}
|v_{w,n}^{\vep}(z)|\leq \frac{\Vert
v_{w,n}^{\vep}\Vert_2}{\sqrt{\pi} (1-|z|^2)}\leq \frac{\Vert
K_1^{\vep}\Vert_{\infty}^{2n-2}\Vert
\mu_2^{\vep}\Vert_{\infty}}{\pi (1-|z|^2)(1-|w|^2)}
\end{align*}
Choose $C_1$ and $C_2$ such that $\Vert K_1 \Vert_{\infty} < C_1
<1$ and $\Vert \mu_2 \Vert_{\infty} <C_2<1$. By the continuity of
the map $\hat{\cP}:T(1)\rightarrow \cB(\ell^2)$ proved in the
Appendix A of \cite{TT3}, the canonical complex analytic embedding
$T(\Gamma_1)\rightarrow T(1)$ and the smooth dependence of $\mu$
on $\vep$, we can find a number $r>0$, such that for all $\vep$ in
the ball $\{ |\vep| < r\}$, we have $\Vert
K_1^{\vep}\Vert_{\infty}<C_1$ and $\Vert
\mu_2^{\vep}\Vert_{\infty} <C_2$. Hence for $|\vep|<r$,
\begin{align*}
\left|\frac{\pa}{\pa \vep} \iint\limits_{\Gamma_1\bk \Del}
\frac{\mathcal{K}_{1,n}(z,z)}{n}d^2z \right|=&
\left|\;\iint\limits_{\Gamma_1\bk \Del}
v_{z,n}^{\vep}(z)d^2z\right| \leq
\frac{C_1^{2n-2}C_2}{\pi}\iint\limits_{\Gamma_1\bk \Del}
\frac{d^2z}{(1-|z|^2)^2}\\
=&\frac{C_1^{2n-2}C_2}{4\pi}\text{Area}(\Gamma_1\bk \Del).
\end{align*}
Consequently, by Weierstrass-M-test, the series
\begin{align*}
\sum_{n=1}^{\infty} \frac{1}{n}
\frac{\pa}{\pa\vep}\iint\limits_{\Gamma_1\bk \Del}
\mathcal{K}_{1,n}^{\vep}(z,z)d^2z
\end{align*}
converges uniformly and absolutely on the set $\{\vep\in \C\;;\;
|\vep|<r\}$. The same proof above shows that as a function of $z$,
the series
\begin{align*}
\sum_{n=1}^{\infty} v_{z,n}^{\vep}(z)
\end{align*}
converges uniformly on any compact subset of $\Del$; in
particular, on a fundamental domain of $\Gamma_1$ on $\Del$.
Therefore,
\begin{align*}
&\sum_{n=1}^{\infty} \frac{1}{n} \iint\limits_{\Gamma_1\bk \Del}
\frac{\pa}{\pa\vep}\mathcal{K}_{1,n}^{\vep}(z,z)d^2z=-\sum_{n=1}^{\infty}
\iint\limits_{\Gamma_1\bk \Del}v_{z,n}^{\vep}(z)
d^2z=-\iint\limits_{\Gamma_1\bk
\Del}\sum_{n=1}^{\infty}v_{z,n}^{\vep}(z)
d^2z\\
=&-\iint\limits_{\Gamma_1\bk
\Del}\iint\limits_{\Del}\iint\limits_{\Del^*}\iint\limits_{\Del^*}\mu_2^{\vep}
(u)K_2^{\vep}(z,u)K_4^{\vep}(u,\eta)\\
&\hspace{3cm}(K_2^{\vep})^*(\eta,\zeta)
\left(\sum_{n=1}^{\infty}\mathcal{K}_{1,
n-1}^{\vep}(\zeta,z)\right)d^2u d^2\eta d^2\zeta d^2z.
\end{align*}
The conclusion of the lemma then follows from the standard
operator theory that $$\sum_{n=1}^{\infty}\mathcal{K}_{1,
n-1}^{\vep}(\zeta,z)=(I-\KKK_1^{\vep})^{-1}(\zeta,z).$$

\end{proof}

Now we state a lemma we need from elementary analysis:
\begin{lemma}\label{lemma3}
Let $\mathcal{O}$ be a ball with center at the origin of $\C$ and
let $h_n : \mathcal{O}\rightarrow \R$ be a sequence of
differentiable real-valued functions on $\mathcal{O}$ that
converges to the function $h: \mathcal{O}\rightarrow \R$. If
$\frac{\pa h}{\pa z}:\mathcal{O}\rightarrow \C$ converges
uniformly to $k: \mathcal{O}\rightarrow \C$, then
$$\frac{\pa h}{\pa z}= k.$$
\end{lemma}
Given $[\mu]\in \mathfrak{D}_g$, let $\vartheta([\mu])\in
\Omega^{2,0}(\Gamma^{\mu})$ be the quadratic differential defined
by
\begin{align*}
\vartheta([\mu]) (z) = \begin{cases} \mathcal{S}(f^{-1})(z),
\hspace{1cm}&\text{if}\; z\in \Omega_1,\\
\mathcal{S}(g^{-1})(z), &\text{if}\;z\in\Omega_2.
\end{cases}
\end{align*}
Here
$$\mathcal{S}(h)=\left(\frac{h^{\prime\prime}}{h^{\prime}}\right)^{\prime}-\frac{1}{2}
\left(\frac{h^{\prime\prime}}{h^{\prime}}\right)^2$$ is the
Schwarzian derivative of the function $h$.
 We have

\begin{theorem}\label{vary}
The real-valued function $S_2: \mathfrak{D}_g\rightarrow \R$ is a
differentiable function. At the point $[\nu]\in \mathfrak{D}_g$,
its variation along the direction $\mu\in
\Omega^{-1,1}(\Gamma^{\nu})$ is given by
\begin{displaymath}
\frac{\pa S_2}{\pa \vep_{ \mu}}([\nu])=-
\frac{1}{6\pi}\iint\limits_{\Gamma^{\nu}\bk\Omega}
\vartheta([\nu])\mu.
\end{displaymath}
\end{theorem}
\begin{proof}

First, we assume $\mu$ has support on $\Omega_2$. From Lemma
\ref{lemma3} and Lemma \ref{lemma4}, we have
\begin{align*}
&\frac{\pa S_2}{\pa\vep_{\mu}}[\nu]= \frac{\pa}{\pa
\vep}\Bigr|_{\vep=0} \iint\limits_{\Gamma_1\bk
\Del}\mathfrak{O}_1^{\vep}(z,z)d^2z\\=&
-\iint\limits_{\Del^*}\iint\limits_{\Gamma_1\bk\Del}
\iint\limits_{\Del} \iint\limits_{\Del^*}
\mu_2(u)K_2(z,u) K_4(u,\eta)\\
&\hspace{5cm} K_2^*(\eta, \zeta) (I-\KKK_1)^{-1}(\zeta, z)d^2\eta
d^2\zeta d^2z d^2u.
\end{align*}
Define $R_2:\Del^*\times \Del\rightarrow\C$ and
$R_2^*:\Del\times\Del^*\rightarrow \C$ be as in the proof of
Theorem 3.1 in \cite{TT3}. Then
\begin{align*}
&\frac{\pa
S_2}{\pa\vep_{\mu}}[\nu]=-\iint\limits_{\Del^*}\iint\limits_{\Gamma_1\bk\Del}
\iint\limits_{\Del} \iint\limits_{\Del^*}
\iint\limits_{\Del^*}\mu_2(u)K_2(z,u) K_4(u,\eta)\\
&\hspace{5cm} K_2^*(\eta, \zeta) R_2^*(\zeta,
v)R_2(v,z)d^2vd^2\eta d^2\zeta d^2z d^2u.
\end{align*} By the transformation property of
the functions $K_l$, $l=1,2,3,4$ with respect to the groups
$\Gamma_1$ and $\Gamma_2$, we can transform the integral into
\begin{align*}
&-\iint\limits_{\Gamma_2\bk\Del^*}\iint\limits_{\Del}
\iint\limits_{\Del} \iint\limits_{\Del^*}
\iint\limits_{\Del^*}\mu_2(u)K_2(z,u) K_4(u,\eta)\\
&\hspace{3cm} K_2^*(\eta, \zeta) R_2^*(\zeta,
v)R_2(v,z)d^2vd^2\eta d^2\zeta d^2z d^2u,
\end{align*}
which can be further be manipulated as in the proof of Theorem 3.1
in \cite{TT3} to get
\begin{displaymath}
\frac{\pa S_2}{\pa \vep_{
\mu}}([\nu])=\frac{1}{6\pi}\iint\limits_{\Gamma_2\bk\Del^*}\mathcal{S}(g)\mu_2=-
\frac{1}{6\pi}\iint\limits_{\Gamma^{\nu}\bk\Omega_2}
\mathcal{S}(g^{-1})\mu.
\end{displaymath}

Now if $\mu$ has support on $\Omega_1$, we let $\kappa=
\iota^*\nu$ and $\eta=\iota^*\mu$. Then $\eta\in
\Omega^{-1,1}(\Gamma^{\kappa})$ and it has support on
$\Omega_2[\kappa]$. By Theorem \ref{invstable} and what we have
proved above
\begin{align*}
\frac{\pa S_2}{\pa \vep_{\mu}}([\nu])=\frac{\pa S_2}{\pa
\vep_{\eta}}([\kappa])=
-\frac{1}{6\pi}\iint\limits_{\Gamma^{\kappa}\bk\Omega_2[\kappa]}
\mathcal{S}(g[\kappa]^{-1})
\eta=-\frac{1}{6\pi}\iint\limits_{\Gamma^{\nu}\bk
\Omega_1}\mathcal{S}(f[\nu]^{-1})\mu.
\end{align*}

Finally, for general $\mu\in \Omega^{-1,1}(\Gamma^{\nu})$, we
write $\mu= \alpha+\beta$, where $\alpha$ has support on
$\Omega_1$ and $\beta$ has support on $\Omega_2$. Then
$$\frac{\pa S_2}{\pa \vep_{\mu}}([\nu])=
\frac{\pa S_2}{\pa \vep_{\alpha}}([\nu])+ \frac{\pa S_2}{\pa
\vep_{\beta}}([\nu]),$$ so the result of the theorem follows.
\end{proof}

\section{The classical Liouville action}
Let $\Gamma'$ be a normalized cocompact Fuchsian group of genus
$g$ realized as a subgroup of $\PSL(2,\R)$
 and let $\alpha(z)= (z-i)/(z+i)$ be the linear fractional
transformation that maps the upper half plane $\U$ to the unit
disc $\Del$. The Fuchsian group $\Gamma= \alpha\circ
\Gamma'\circ\alpha^{-1}$ is then a subgroup of $\PSU(1,1)$.
$\alpha$ induces a complex analytic isomorphism $\mathfrak{G}:
\mathfrak{D}(\Gamma')\rightarrow \mathfrak{D}(\Gamma)$ by
$$[\mu]\in \mathfrak{D}(\Gamma')\mapsto \mathfrak{G}([\mu])= [
(\alpha^{-1})^*\mu].$$ We define the function $\SSS_2:
\mathfrak{D}(\Gamma')\rightarrow \R$ by $\SSS_2=S_2\circ
\mathfrak{G}$.

Given $[\nu]\in\mathfrak{D}(\Gamma')$, we let $w_{\nu}$ be the
unique q.c. mapping with Beltrami differential $\nu$ and fixing
the points $0,1,\infty$. Let $\kappa= (\alpha^{-1})^*\nu$, then
$w_{\nu}= \alpha^{-1} \circ w_{\kappa}\circ \alpha$. Let $J_1=
\alpha^{-1}\circ f\circ \alpha$, $J_2= \alpha^{-1} \circ g \circ
\alpha$, where $(f,g)$ is the pair of univalent functions
associated to $[\kappa]\in \mathfrak{D}(\Gamma)$.

\subsection{Relation between classical Liouville action and the function $\SSS_2$}
 Given $[\mu]\in
\mathfrak{D}(\Gamma')$, let $\tilde{\vartheta}([\mu])\in
\Omega^{2,0}((\Gamma')^{\mu})$ be the quadratic differential
defined by
\begin{align*}
\tilde{\vartheta}([\mu]) (z) = \begin{cases}
\mathcal{S}(J_1^{-1})(z),
\hspace{0.5cm}&\text{if}\;\; z\in w_{\mu}(\U),\\
\mathcal{S}(J_2^{-1})(z), &\text{if}\;\;z\in w_{\mu}(\bL).
\end{cases}
\end{align*}
Here $\bL$ is the lower half plane. If $\kappa=
(\alpha^{-1})^*\mu$, then $\vartheta([\kappa])=
(\alpha^{-1})^*\tilde{\vartheta}([\mu])$. It follows from Theorem
\ref{vary} that
\begin{theorem}\label{vary2}
Given a subgroup $\Gamma'$ of $\PSL(2,\R)$, which is a normalized
cocompact Fuchsian group of genus $g$, the real-valued function
$\SSS_2: \mathfrak{D}(\Gamma')\rightarrow \R$ is a differentiable
function. At the point $[\nu]\in \mathfrak{D}(\Gamma')$, its
variation along the direction $\mu\in
\Omega^{-1,1}((\Gamma')^{\nu})$ is given by
\begin{displaymath}
\frac{\pa \SSS_2}{\pa \vep_{ \mu}}([\nu])=-
\frac{1}{6\pi}\iint\limits_{(\Gamma')^{\nu}\bk\Omega'}
\tilde{\vartheta}([\nu])\mu.
\end{displaymath}
Here $\Omega' =w_{\nu}(\U)\sqcup w_{\nu}(\bL)$ is the set of
discontinuity of the group $(\Gamma')^{\nu}$.
\end{theorem}
\begin{proof}
We let $\kappa = (\alpha^{-1})^*\nu$, $\eta= (\alpha^{-1})^*\mu$.
Then by Theorem \ref{vary},
\begin{align*}
\frac{\pa \SSS_2}{\pa \vep_{ \mu}}([\nu])= \frac{\pa S_2}{\pa
\vep_{ \eta}}([\kappa])=-
\frac{1}{6\pi}\iint\limits_{\Gamma^{\kappa}\bk\Omega}
\vartheta([\kappa])\eta=-
\frac{1}{6\pi}\iint\limits_{(\Gamma')^{\nu}\bk\Omega'}
\tilde{\vartheta}([\nu])\mu.
\end{align*}
\end{proof}

In \cite{TT1}, we define the classical Liouville action $S_{cl}:
\mathfrak{D}(\Gamma')\rightarrow \R$ and prove that $-S_{cl}$ is a
Weil-Petersson potential on $\mathfrak{D}(\Gamma')$.  Theorem 4.2
in \cite{TT1} says that
\begin{theorem}\label{vary3}
Let $\Gamma'$ be as in Theorem \ref{vary2}. At the point $[\nu]\in
\mathfrak{D}(\Gamma')$, the variation of the function $S_{cl}$
along the direction $\mu\in \Omega^{-1,1}((\Gamma')^{\nu})$ is
given by
\begin{displaymath}
\frac{\pa S_{cl}}{\pa \vep_{ \mu}}([\nu])=
2\iint\limits_{(\Gamma')^{\nu}\bk\Omega'}
\tilde{\vartheta}([\nu])\mu.
\end{displaymath}

\end{theorem}

Theorem \ref{vary2} and Theorem \ref{vary3} give a relation
between $S_{cl}$ and $\SSS_2$.
\begin{theorem}
Let $\Gamma'$ be as in Theorem \ref{vary2}. On the deformation
space $\mathfrak{D}(\Gamma')$, we have
\begin{align*}
S_{cl} = -12\pi \SSS_2+8\pi (2g-2).
\end{align*}
\end{theorem}
\begin{proof}
Since $\mathfrak{D}(\Gamma')$ is connected, from Theorem
\ref{vary2} and Theorem \ref{vary3}, we have
$$S_{cl} = -12\pi \SSS_2 +C,$$ where $C$ is a constant. Now at the
origin $[0]$ of $\mathfrak{D}(\Gamma')$, $S_{cl}([0])=8\pi(2g-2)$
and $\SSS_2([0])=0$, hence $C=8\pi(2g-2)$.
\end{proof}

This gives us the following inequality for the classical Liouville
action.
\begin{corollary}\label{inequality}
The classical Liouville action
$S_{cl}:\mathfrak{D}(\Gamma')\rightarrow \R$ satisfies the
following inequality
$$ S_{cl}\leq 8\pi (2g-2).$$
It attains its maximum value along the subspace of Fuchsian
groups.
\end{corollary}
\begin{proof}
This follows from the theorem above and the nonnegativity of $S_2$
established in Corollary \ref{pos}.
\end{proof}
\begin{remark}
It follows that the normalized potential $-S_{cl}+8\pi(2g-2)$ on
the quasi-Fuchsian deformation space $\mathfrak{D}(\Gamma')$ is a
nonnegative function.
\end{remark}

\bibliographystyle{amsalpha}
\bibliography{positive}
\end{document}